\documentclass[11pt,a4paper,english]{article}
\usepackage{mathptmx}
\usepackage[T1]{fontenc}
\usepackage[latin9]{inputenc}
\usepackage{color}
\usepackage{array}
\usepackage{amsmath}
\usepackage{amssymb}
\usepackage{graphicx}
\usepackage{esint}
\usepackage{listings}
\makeatletter

\pdfpageheight\paperheight
\pdfpagewidth\paperwidth

\usepackage{url}
\DeclareUrlCommand\email{}

\usepackage{iamPreprintStyle}

\makeatother

\usepackage{babel}
\begin{document}
\preprintTitle{A Matlab Tutorial for Diffusion-Convection-Reaction Equations using DGFEM} \preprintAuthor{Murat Uzunca%
\footnote{Department of Mathematics, Middle East Technical University, 06800 Ankara, Turkey, \textit{email}: uzunca@metu.edu.tr%
}, B\"{u}lent Karas\"{o}zen%
\footnote{Department of Mathematics and Institute of Applied Mathematics, Middle East Technical University, 06800 Ankara, Turkey, \textit{email}: bulent@metu.edu.tr%
}} 
\preprintAbstract{\small We present a collection of MATLAB routines using discontinuous Galerkin finite elements method (DGFEM) for solving steady-state diffusion-convection-reaction equations. The code employs the sparse matrix facilities of MATLAB with "vectorization" and uses multiple matrix multiplications {\it "MULTIPROD"} \cite{leva} to increase the efficiency of the program. } 

\preprintKeywords{Discontinuous Galerkin FEMs, Diffusion-convection-reaction equations, Matlab} 

\preprintDate{August 2014} 
\preprintNo{2014-4} 
\makePreprintCoverPage


\section{DG discretization of the linear model problem}

Many engineering problems such as chemical reaction processes, heat conduction,
 nuclear reactors, population dynamics etc. are governed by convection-diffusion-reaction partial
differential equations (PDEs). The general model problem used in the code is
\begin{subequations}\label{1}
\begin{align}
  \alpha u - \epsilon\Delta u + {\bf b}\cdot\nabla u  &= f \quad \; \text{ in } \; \Omega, \\
   u &= g^D  \quad \text{on } \; \Gamma^D, \\
	\epsilon\nabla u\cdot {\bf n} &= g^N  \quad \text{on } \; \Gamma^N.
\end{align}
\end{subequations}
The domain $\Omega$ is bounded, open, convex in $\mathbb{R}^2$ with boundary $\partial\Omega =\Gamma^D\cup\Gamma^N$, $\Gamma^D\cap\Gamma^N=\emptyset$, $0<\epsilon\ll 1 $ is the diffusivity constant, $f\in L^2(\Omega)$ is the source function, ${\bf b}\in\left(W^{1,\infty}(\Omega)\right)^2$ is the velocity field, $g^D\in H^{3/2}(\Gamma^D )$ is the Dirichlet boundary condition, $g^N\in H^{1/2}(\Gamma^N )$ is the Neumann boundary condition and ${\bf n}$ denote the unit outward normal vector to the boundary.\\

\noindent The weak formulation of  (\ref{1}) reads as: find $u\in U$ such that

\begin{equation} \label{2}
\int_{\Omega}(\epsilon\nabla u\cdot\nabla v+{\bf b}\cdot\nabla uv+\alpha uv)dx = \int_{\Omega}fvdx +\int_{\Gamma^N}g^Nv ds\; , \quad \forall v\in V
\end{equation}
where the solution space $U$ and the test function space $V$ are given by

\begin{equation*}
U= \{ u\in H^1(\Omega) \, : \; u=g^D \text{ on } \Gamma^D\}, \quad
V= \{ v\in H^1(\Omega) \, : \; v=0 \text{ on } \Gamma^D\}.
\end{equation*}
Being the next step, an approximation to the problem (\ref{2}) is found in a finite-dimensional space $V_h$. In case of classical (continuous) FEMs, the space $V_h$ is set to be the set of piecewise continuous polynomials vanishing on the boundary $\partial\Omega$. In contrast to the continuous FEMs, the DGFEMs uses the set of piecewise polynomials that are fully discontinuous at the interfaces. In this way,  the DGFEMs approximation allows to capture the sharp gradients or singularities that affect the numerical solution locally. Moreover, the functions in $V_h$ do not need to vanish at the boundary since the boundary conditions in DGFEMs are imposed weakly.\\

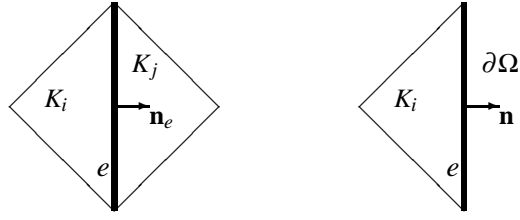
\begin{figure}
\centering
\setlength{\unitlength}{2.3mm}
\begin{picture}(35, 13)
\put(5,7){\line(1,1){6}}
\put(5,7){\line(1,-1){6}}
\put(17,7){\line(-1,1){6}}
\put(17,7){\line(-1,-1){6}}
\put(7,7){\small$K_{i}$ }
\put(12,9){\small$K_{j}$ }
\put(10,3){$e$ }
\put(11,7){\vector(1,0){2}}
\put(13,6){\small${\bf n}_{e}$}

\put(25,7){\line(1,1){6}}
\put(25,7){\line(1,-1){6}}
\put(27,7){\small$K_{i}$ }
\put(32,9){\small$\partial\Omega$ }
\put(30,3){$e$ }
\put(31,7){\vector(1,0){2}}
\put(33,6){\small${\bf n}$}

\linethickness{0.7mm} \put(11,1){\line(0,1){12}}
\linethickness{0.7mm} \put(31,1){\line(0,1){12}}
\end{picture}
\caption{Two adjacent elements sharing an edge (left); an element near to domain boundary (right)}
\label{jump}
\end{figure}

\noindent In our code, the discretization of the problem (\ref{1}) is based on the discontinuous Galerkin methods for the diffusion part \cite{arnold02uad,riviere08dgm} and the upwinding  for the convection part \cite{ayuso09dgm,houston02dhp}. Let $\{\xi_h\}$ be a family of shape regular meshes with the elements (triangles) $K_i\in\xi_h$ satisfying $\overline{\Omega}=\cup \overline{K}$ and $K_i\cap K_j=\emptyset$ for $K_i$, $K_j$ $\in\xi_h$. Let us denote by $\Gamma_0$, $\Gamma_D$ and $\Gamma_N$ the set of interior, Dirichlet boundary and Neumann boundary edges, respectively, so that $\Gamma_0\cup\Gamma_D\cup\Gamma_N$ forms the skeleton of the mesh. For any $K\in\xi_h$, let $\mathbb{P}_k(K)$ be the set of all polynomials of degree at most $k$ on $K$. Then, set the finite dimensional solution and test function space by
$$
V_h=\left\{ v\in L^2(\Omega ) : v|_{K}\in\mathbb{P}_k(K) ,\; \forall K\in \xi_h \right\}\not\subset V.
$$
Note that the trial and test function spaces are the same because the boundary conditions in discontinuous Galerkin methods are imposed in a weak manner. Since the functions in $V_h$ may have  discontinuities along the inter-element boundaries, along an interior edge, there would be two different traces from the adjacent elements sharing that edge. In the light of this fact, let us first introduce some notations before giving the DG formulation. Let $K_i$, $K_j\in\xi_h$ ($i<j$) be two adjacent elements sharing an interior edge $e=K_i\cap K_j\subset \Gamma_0$ (see Fig.\ref{jump}). Denote the trace of a scalar function $v$ from inside $K_i$ by $v_{i}$ and from inside $K_j$ by $v_{j}$. Then, set the jump and average values of $v$ on the edge $e$
$$
[v]= v_{i}{\bf n}_e- v_{j}{\bf n}_e , \quad \{ v\}=\frac{1}{2}(v_{i}+ v_{j}),
$$
where ${\bf n}_e$ is the unit normal to the edge $e$ oriented from $K_i$ to $K_j$. Similarly, we set the jump and average values of a vector valued function ${\bf q}$ on e
$$
[{\bf q}]= {\bf q}_{i}\cdot {\bf n}_e- {\bf q}_{j}\cdot {\bf n}_e , \quad \{ {\bf q}\}=\frac{1}{2}({\bf q}_{i}+ {\bf q}_{j}),
$$
Observe that $[v]$ is a vector for a scalar function $v$, while, $[{\bf q}]$ is scalar for a vector valued function ${\bf q}$. On the other hands, along any boundary edge $e=K_i\cap \partial\Omega$, we set
$$
[v]= v_{i}{\bf n} , \quad \{ v\}=v_{i}, \quad [{\bf q}]={\bf q}_{i}\cdot {\bf n}, \quad \{ {\bf q}\}={\bf q}_{i}
$$
where ${\bf n}$ is the unit outward normal to the boundary at $e$.

\noindent We also introduce the inflow parts of the domain boundary and the boundary of a mesh element $K$, respectively
$$
\Gamma^- = \{ x\in\partial \Omega \; : \; {\bf b}(x)\cdot {\bf n}(x)<0\} \; , \quad \partial K^- = \{ x\in\partial K \; : \; {\bf b}(x)\cdot {\bf n}_K(x)<0\}.
$$ 
Then, the DG discretized system to the problem (\ref{1}) combining with the upwind discretization for the convection part reads as: find $u_h\in V_h$ such that
\begin{equation} \label{ds}
a_{h}(u_{h},v_{h})=l_{h}(v_{h}) \qquad \forall v_h\in V_h,
\end{equation}
\begin{subequations}
\begin{align}
a_{h}(u_{h}, v_{h})=& \sum \limits_{K \in {\xi}_{h}} \int_{K} \epsilon \nabla u_{h}\cdot\nabla v_{h} dx + \sum \limits_{K \in {\xi}_{h}} \int_{K} ({\bf b} \cdot \nabla u_{h}+\alpha u_h) v_{h} dx\nonumber \\
&- \sum \limits_{ e \in \Gamma_{0}\cup\Gamma_{D}}\int_{e} \{\epsilon \nabla u_{h}\} \cdot [v_{h}] ds  + \kappa \sum \limits_{ e \in \Gamma_{0}\cup\Gamma_{D}} \int_{e} \{\epsilon \nabla v_{h}\} \cdot [u_{h}] ds \nonumber \\
&+ \sum \limits_{K \in {\xi}_{h}}\int_{\partial K^-\setminus\partial\Omega } {\bf b}\cdot {\bf n} (u_{h}^{out}-u_{h}^{in})  v_{h} ds - \sum \limits_{K \in {\xi}_{h}} \int_{\partial K^-\cap \Gamma^{-}} {\bf b}\cdot {\bf n} u_{h}^{in} v_{h}  ds \nonumber  \\
&+ \sum \limits_{ e \in \Gamma_{0}\cup\Gamma_{D}}\frac{\sigma \epsilon}{h_{e}} \int_{e} [u_{h}] \cdot [v_{h}] ds, \nonumber  \\
l_{h}( v_{h})=&  \sum \limits_{K \in {\xi}_{h}} \int_{K} f v_{h} dx
+ \sum \limits_{e \in \Gamma_{D}} \int_e g^D \left( \frac{\sigma \epsilon}{h_{e}} v_{h} -  {\epsilon\nabla v_{h}} \cdot {\bf n} \right) ds  \nonumber \\
&- \sum \limits_{K \in {\xi}_{h}}\int_{\partial K^-\cap \Gamma^{-}} {\bf b}\cdot {\bf n} g^D v_{h}  ds + \sum \limits_{e \in \Gamma_{N}} \int_e g^N v_{h} ds, \nonumber
\end{align}
\end{subequations}
where $u_{h}^{out}$ and $u_{h}^{in}$ denotes the values on an edge from outside and inside of an element $K$, respectively. The parameter $\kappa$ determines the type of DG method, which takes the values $\{ -1,1,0 \}$: $\kappa =-1$ gives {\it "symmetric interior penalty Galerkin"} (SIPG) method, $\kappa =1$ gives {\it "non-symmetric interior penalty Galerkin"} (NIPG) method and $\kappa =0$ gives {\it "inconsistent interior penalty Galerkin"} (IIPG) method. The parameter $\sigma\in\mathbb{R}_0^+$ is called the penalty parameter which should be sufficiently large; independent of the mesh size $h$ and the diffusion coefficient $\epsilon$ \cite{riviere08dgm} [Sec. 2.7.1]. In our code, we choose the penalty parameter $\sigma$ on interior edges depending on the polynomial degree $k$ as $\sigma=3k(k+1)$ for the SIPG and IIPG methods, whereas,  we take $\sigma=1$ for the NIPG method. On boundary edges, we take the penalty parameter as twice of the penalty parameter on interior edges.

\section{Descriptions of the MATLAB code}

The given codes are mostly self-explanatory with comments to explain what each section of the code does. In this section, we give the description of our main code. The use of the main code consists of three parts
\begin{enumerate}
\item Mesh generation,
\item Entry of user defined quantities (boundary conditions, order of basis etc.),
\item Forming and solving the linear systems,
\item Plotting the solutions.
\end{enumerate}

\noindent Except the last one, all the parts above take place in the m-file {\it Main\_Linear.m} which is the main code to be used by the users for linear problems without need to entry to any other m-file. The last part, plotting the solutions, takes place in the m-file {\it dg\_error.m}.

\subsection{Mesh generation}

In this section, we define the data structure of a triangular mesh on a polygonal domain
in $\mathbb{R}^2$. The data structure presented here is based on simple arrays \cite{chen08fem} which are stored in a MATLAB "struct" that collects two or more data fields in
one object that can then be passed to routines. To obtain an initial mesh, firstly, we define the nodes, elements, Dirichlet and Neumann conditions in the m-file {\it Main\_Linear.m}, and we call the {\it getmesh} function to form the initial mesh structure {\it mesh}.

\begin{lstlisting}
% Generate the mesh
% Nodes
Nodes = [0,0;0.5,0;1,0;0,0.5;0.5,0.5;1,0.5;0,1;0.5,1;1,1];
% Elements
Elements = [4,1,5;1,2,5;5,2,6; 2,3,6;7,4,8;4,5,8;8,5,9;5,6,9];
% Dirichlet bdry edges
Dirichlet = [1,2;2,3;1,4;3,6;4,7;6,9;7,8;8,9];
% Neumann bdry edges
Neumann = [];
% Initial mesh struct
mesh = getmesh(Nodes,Elements,Dirichlet,Neumann);
\end{lstlisting}

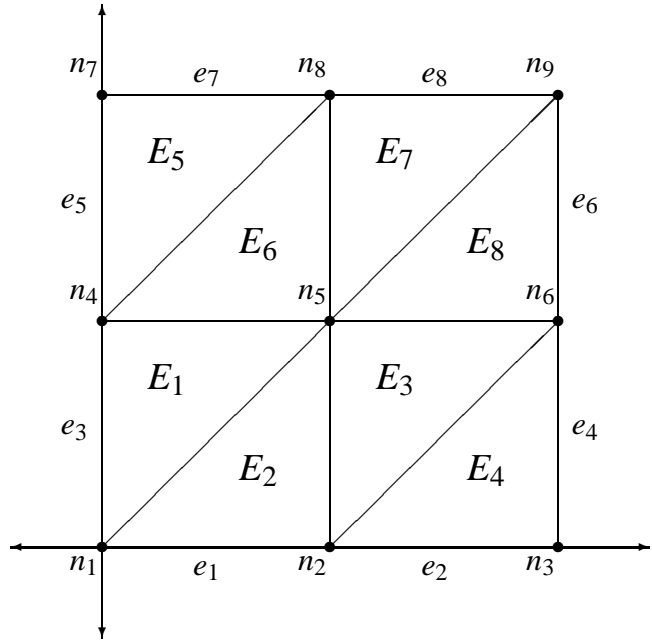
\begin{figure}
\centering
\setlength{\unitlength}{6mm}
\begin{picture}(14, 14)
\put(2,2){\vector(1,0){12}}
\put(2,2){\vector(0,1){12}}
\put(2,2){\vector(-1,0){2}}
\put(2,2){\vector(0,-1){2}}
\put(2,12){\line(1,0){10}}
\put(2,7){\line(1,0){10}}
\put(7,2){\line(0,1){10}}
\put(12,2){\line(0,1){10}}

\put(2,2){\line(1,1){10}}
\put(2,7){\line(1,1){5}}
\put(7,2){\line(1,1){5}}

\put(2,2){\circle*{0.25}}
\put(2,7){\circle*{0.25}}
\put(2,12){\circle*{0.25}}
\put(7,2){\circle*{0.25}}
\put(7,7){\circle*{0.25}}
\put(7,12){\circle*{0.25}}
\put(12,2){\circle*{0.25}}
\put(12,7){\circle*{0.25}}
\put(12,12){\circle*{0.25}}

\put(1.3,1.5){\large$n_1$ }
\put(6.3,1.5){\large$n_2$ }
\put(11.3,1.5){\large$n_3$ }
\put(1.3,7.5){\large$n_4$ }
\put(6.3,7.5){\large$n_5$ }
\put(11.3,7.5){\large$n_6$ }
\put(1.3,12.5){\large$n_7$ }
\put(6.3,12.5){\large$n_8$ }
\put(11.3,12.5){\large$n_9$ }

\put(3,5.5){\Large$E_1$ }
\put(5,3.5){\Large$E_2$ }
\put(8,5.5){\Large$E_3$ }
\put(10,3.5){\Large$E_4$ }
\put(3,10.5){\Large$E_5$ }
\put(5,8.5){\Large$E_6$ }
\put(8,10.5){\Large$E_7$ }
\put(10,8.5){\Large$E_8$ }

\put(4,1.4){\large$e_1$ }
\put(9,1.4){\large$e_2$ }
\put(1.1,4.5){\large$e_3$ }
\put(12.3,4.5){\large$e_4$ }
\put(1.1,9.5){\large$e_5$ }
\put(12.3,9.5){\large$e_6$ }
\put(4,12.3){\large$e_7$ }
\put(9,12.3){\large$e_8$ }

\end{picture}
\caption{Initial mesh on the unit square $\Omega = [0,1]^2$ with nodes $n_i$, triangles $E_j$ and edges $e_k$}
\label{meshdata}
\end{figure}

\noindent As it can be understood, each row in the {\bf Nodes} array corresponds to a mesh node with the first column keeps the x-coordinate of the node and the second is for the y-coordinate, and the $i-th$ row of the {\bf Nodes} array is called the node having index $i$. In the {\bf Elements} array, each row with 3 columns corresponds to a triangular element in the mesh containing the indices of the nodes forming the 3 vertices of the triangles in the counter-clockwise orientation. Finally, in the {\bf Dirichlet} and {\bf Neumann} arrays, each row with 2 columns corresponds to a Dirichlet and Neumann boundary edge containing the indices of the starting and ending nodes, respectively (see Fig.\ref{meshdata}).
The mesh "struct" in the code has the following fields:
\begin{itemize}
	\item Nodes, Elements, Edges, intEdges, DbdEdges, NbdEdges, intEdges
	\item vertices1, vertices2, vertices3,
	\item Dirichlet, Neumann, EdgeEls, ElementsE.
\end{itemize}
which can be reached by {\it mesh.Nodes}, {\it mesh.Elements} and so on, and they are used by the other functions to form the DG construction. The initial mesh can be uniformly refined several times in a "for loop" by calling the function {\it uniformrefine}.

\begin{lstlisting}
for jj=1:2
 mesh=uniformrefine(mesh); %Refine mesh
end
\end{lstlisting}

\subsection{User defined quantities}

There are certain input values that have to be supplied by the user. Here, we will describe that how one can define these quantities in the main code {\it Main\_Linear.m}. \\

\noindent One determines the type of the DG method (SIPG, NIPG or IIPG) and the order of the polynomial basis to be used by the variables {\it method} and {\it degree}, respectively. According to these choices, the values of the penalty parameter and the parameter $\kappa \in \{-1,1,0\}$ defining DG method in (\ref{ds}) are set by calling the sub-function {\it set\_parameter}.

\begin{lstlisting}
% method : NIPG=1, SIPG=2, IIPG=3
method=2;

% Degree of polynomials
degree=1;

% Set up the problem
[penalty,kappa]=set_parameter(method,degree);
\end{lstlisting}
The next step is to supply the problem parameters. The diffusion constant $\epsilon$, the advection vector $\bf b$ and the linear reaction term $\alpha$ are defined via the sub-functions {\it fdiff}, {\it fadv} and {\it freact}, respectively.

\begin{lstlisting}
%% Define diffusion, advection, and reaction as subfunctions

% Diffusion
function diff = fdiff(x,y)
  diff = (10^(-6)).*ones(size(x));
end

% Advection
function [adv1,adv2] = fadv(x,y)
  adv1 =(1/sqrt(5))*ones(size(x));
  adv2 =(2/sqrt(5))*ones(size(x));
end

% Linear reaction
function react = freact(x,y)
  react = ones(size(x));
end
\end{lstlisting}
The exact solution (if exists) and the source function $f$ are defined via the sub-functions {\it fexact} and {\it fsource}, respectively. Finally, the boundary conditions are supplied via the sub-functions {\it DBCexact} and {\it NBCexact}. \\

\begin{lstlisting}
% First derivative wrt x
yex_x=(-1./(sqrt(5*diff))).*(sech((2*x-y-0.25)./...
    (sqrt(5*diff)))).^2;
% First derivative wrt y
yex_y=((0.5)./(sqrt(5*diff))).*(sech((2*x-y-0.25)./...
    (sqrt(5*diff)))).^2;
% Second derivative wrt x
yex_xx=((0.8)./diff).*tanh((2*x-y-0.25)./(sqrt(5*diff))).*...
    (sech((2*x-y-0.25)./(sqrt(5*diff)))).^2;
% Second derivative wrt y
yex_yy=((0.2)./diff).*tanh((2*x-y-0.25)./(sqrt(5*diff))).*...
    (sech((2*x-y-0.25)./(sqrt(5*diff)))).^2;
		
% Force function
source=-diff.*(yex_xx+yex_yy)+(adv1.*yex_x+adv2.*yex_y)+...
\end{lstlisting}

\subsection{Forming and solving linear systems}

To form the linear systems, firstly, let us rewrite the discrete DG scheme (\ref{ds}) as
\begin{equation}\label{dgscheme}
a_{h}(u_{h},v_{h}):= D_{h}(u_{h},v_{h}) + C_{h}(u_{h},v_{h}) + R_{h}(u_{h},v_{h}) =l_{h}(v_{h}) \qquad \forall v_h\in V_h,
\end{equation}
where the forms $D_{h}(u_{h},v_{h})$, $C_{h}(u_{h},v_{h})$ and $R_{h}(u_{h},v_{h})$ corresponding to the diffusion, convection and linear reaction parts of the problem, respectively
\begin{subequations}
\begin{align}
D_{h}(u_{h}, v_{h})=& \sum \limits_{K \in {\xi}_{h}} \int_{K} \epsilon \nabla u_{h}\cdot\nabla v_{h} dx + \sum \limits_{ e \in \Gamma_{0}\cup\Gamma_{D}}\frac{\sigma \epsilon}{h_{e}} \int_{e} [u_{h}] \cdot [v_{h}] ds  \nonumber \\ 
&- \sum \limits_{ e \in \Gamma_{0}\cup\Gamma_{D}}\int_{e} \{\epsilon \nabla u_{h}\} \cdot [v_{h}] ds  + \kappa \sum \limits_{ e \in \Gamma_{0}\cup\Gamma_{D}} \int_{e} \{\epsilon \nabla v_{h}\} \cdot [u_{h}] ds \nonumber \\
C_{h}(u_{h}, v_{h})=&  \sum \limits_{K \in {\xi}_{h}} \int_{K} {\bf b} \cdot \nabla u_{h} v_{h} dx\nonumber \\
&+ \sum \limits_{K \in {\xi}_{h}}\int_{\partial K^-\setminus\partial\Omega } {\bf b}\cdot {\bf n} (u_{h}^{out}-u_{h}^{in})  v_{h} ds - \sum \limits_{K \in {\xi}_{h}} \int_{\partial K^-\cap \Gamma^{-}} {\bf b}\cdot {\bf n} u_{h}^{in} v_{h}  ds \nonumber \\
R_{h}(u_{h}, v_{h})=&  \sum \limits_{K \in {\xi}_{h}} \int_{K} \alpha u_h v_{h} dx \nonumber \\
l_{h}( v_{h})=&  \sum \limits_{K \in {\xi}_{h}} \int_{K} f v_{h} dx
+ \sum \limits_{e \in \Gamma_{D}} \int_e g^D \left( \frac{\sigma \epsilon}{h_{e}} v_{h} -  {\epsilon\nabla v_{h}} \cdot {\bf n} \right) ds  \nonumber \\
&- \sum \limits_{K \in {\xi}_{h}}\int_{\partial K^-\cap \Gamma^{-}} {\bf b}\cdot {\bf n} g^D v_{h}  ds + \sum \limits_{e \in \Gamma_{N}} \int_e g^N v_{h} ds, \nonumber
\end{align}
\end{subequations}
For a set of basis functions $\{\phi_i\}_{i=1}^N$ spanning the space $V_h$, the discrete solution $u_h\in V_h$ is of the form
\begin{equation}\label{sol}
u_h = \sum_{j=1}^N \upsilon_j\phi_j
\end{equation}
where $\upsilon=(\upsilon_1,\upsilon_2, \ldots , \upsilon_N)^T$ is the unknown coefficient vector.
After substituting (\ref{sol}) into (\ref{dgscheme}) and taking $v_h=\phi_i$, we get the linear system of equations
\begin{equation}\label{system}
\sum_{j=1}^N \upsilon_j D_{h}(\phi_j,\phi_i) + \sum_{j=1}^N \upsilon_j C_{h}(\phi_j,\phi_i) + \sum_{j=1}^N \upsilon_j R_{h}(\phi_j,\phi_i) = l_h(\phi_i) \; , \quad i=1,2,\ldots , N
\end{equation}
Thus, for $i=1,2,\ldots , N$, to form the linear system in matrix-vector form, we need the matrices $D, C, R\in \mathbb{R}^{N\times N}$ related to the terms including the forms $D_{h}$, $C_{h}$ and $R_{h}$ in (\ref{system}), respectively, satisfying
$$
D\upsilon + C\upsilon + R\upsilon = F
$$
with the unknown coefficient vector $\upsilon$ and the vector $F\in\mathbb{R}^N$ related to the linear rhs functionals $l_h(\phi_i)$ such that $F_i=l_h(\phi_i)$, $i=1,2,\ldots , N$. In the code {\it Main\_Linear.m}, all the matrices $D,C,R$ and the vector $F$ are obtained by calling the function {\it global\_system}, in which the sub-functions introduced in the previous subsection are used. We set the stiffness matrix, {\it Stiff}, as the sum of the obtained matrices and we solve the linear system for the unknown coefficient vector {\it coef}$:=\upsilon$.

\begin{lstlisting}
%Compute global matrices and rhs global vector
[D,C,R,F]=global_system(mesh,@fdiff,@fadv,@freact,...
      @fsource,@DBCexact,@NBCexact,penalty,kappa,degree);
		
Stiff=D+C+R; % Stiffness matrix

coef=Stiff\F; % Solve the linear system
\end{lstlisting}

\subsection{Plotting the solution}

After solving the problem for the unknown coefficient vector, the solutions are plotted via the the function {\it dg\_error}, and also the $L^2$-error between the exact and numerical solution is computed.

\begin{lstlisting}
% Compute L2-error and plot the solution
[l2err,hmax]=dg_error(coef,mesh,@fexact,@fdiff,degree);
\end{lstlisting}

\section{Models with non-linear reaction mechanisms}

Most of the problems include non-linear source or sink terms. The general model problem in this case is 
\begin{subequations}\label{nonlin}
\begin{align}
  \alpha u - \epsilon\Delta u + {\bf b}\cdot\nabla u + r(u) &= f \quad \; \text{ in } \; \Omega, \\
   u &= g^D  \quad \text{on } \; \Gamma^D, \\
	\epsilon\nabla u\cdot {\bf n} &= g^N  \quad \text{on } \; \Gamma^N.
\end{align}
\end{subequations}
which arises from the time discretization of the time-dependent non-linear diffusion-convection-reaction equations. Here, the coefficient of the linear reaction term, $\alpha >0$, stand for the temporal discretization, corresponding to $1/\Delta t$, where $\Delta t$ is the discrete time-step. The model (\ref{nonlin}) differs from the model (\ref{1}) by the additional non-linear term $r(u)$. To have a unique solution, in addition to the assumptions given in Section 1, we assume that the non-linear reaction term, $r(u)$, is bounded, locally Lipschitz continuous and monotone, i.e. satisfies for any $s, s_1, s_2\ge 0$, $s,s_1, s_2 \in \mathbb{R}$ the following conditions  \cite{uzunca14adg}

\begin{align*}
|r_i(s)| &\leq C , \quad C>0 \\
\| r_i(s_1)-r_i(s_2)\|_{L^2(\Omega)} &\leq L\| s_1-s_2\|_{L^2(\Omega)} , \quad L>0  \\
 r_i\in C^1(\mathbb{R}_0^+), \quad r_i(0) =0, &\quad r_i'(s)\ge 0. 
\end{align*}
The non-linear reaction term $r(u)$ occur in chemical engineering usually in the form of products and rational functions of concentrations, or exponential functions of the temperature, expressed by the Arrhenius law. Such models describe chemical processes and they are strongly coupled as an inaccuracy in one unknown affects all the others.\\

\noindent To solve the non-linear problems, we use the m-file {\it Main\_Nonlinear} which is similar to the m-file {\it Main\_Linear}, but now we use Newton iteration to solve for $i=1,2,\ldots , N$ the non-linear system of equations

\begin{equation}\label{nonlinsystem}
\sum_{j=1}^N \upsilon_j D_{h}(\phi_j,\phi_i) + \sum_{j=1}^N \upsilon_j C_{h}(\phi_j,\phi_i) + \sum_{j=1}^N \upsilon_j R_{h}(\phi_j,\phi_i) + \int_{\Omega} r(u_h)\phi_i dx = l_h(\phi_i) 
\end{equation}
Similar to the linear case, the above system leads to the matrix-vector form 
$$
D\upsilon + C\upsilon + R\upsilon + H(\upsilon) = F
$$
where, in addition to the matrices $D,C,R\in\mathbb{R}^{N\times N}$ and the vector $F\in\mathbb{R}^N$, we also need the vector $H\in\mathbb{R}^N$ related to the non-linear term such that 
$$
H_i(\upsilon) = \int_{\Omega} r\left( \sum_{j=1}^N \upsilon_j \phi_j \right)\phi_i dx \; , \quad i=1,2,\ldots , N.
$$
We solve the nonlinear system by Newton method. For an initial guess $\upsilon^0=(\upsilon^0_1,\upsilon^0_2, \ldots , \upsilon^0_N)^T$, we solve the system
\begin{eqnarray} \label{newton}
J^kw^k &=& -Res^k    \\
 \upsilon^{k+1} &=& w^k + \upsilon^k \; , \quad k=0,1,2,\ldots \nonumber
\end{eqnarray}
until a user defined tolerance is satisfied. In (\ref{newton}), $Res^k$ and $J^k$ denote the vector of system residual and its Jacobian matrix at the current iterate $\upsilon^k$, respectively, given by
\begin{eqnarray*} 
Res^k &=& (D+C+R)\upsilon^k + H(\upsilon^k) - F    \\
 J^k &=& D+C+R + HJ(\upsilon^k)
\end{eqnarray*}
where $HJ(\upsilon^k)$ is the Jacobian matrix of the non-linear vector $H$ at $\upsilon^k$
$$
HJ(\upsilon^k)=
\begin{bmatrix}
\frac{\partial H_1(\upsilon^k)}{\partial \upsilon^k_1} &  \frac{\partial H_1(\upsilon^k)}{\partial \upsilon^k_2} & \cdots & \frac{\partial H_1(\upsilon^k)}{\partial \upsilon^k_N} \\
\vdots & \ddots & & \vdots \\
\frac{\partial H_N(\upsilon^k)}{\partial \upsilon^k_1} &  \frac{\partial H_N(\upsilon^k)}{\partial \upsilon^k_2} & \cdots & \frac{\partial H_N(\upsilon^k)}{\partial \upsilon^k_N}
\end{bmatrix}
$$
In the code {\it Main\_Nonlinear}, obtaining the matrices $D,C,R$ and the rhs vector $F$ is similar to the linear case, but now, additionally, we introduce an initial guess for Newton iteration, and we solve the nonlinear system by Newton method. 

\pagebreak

\begin{lstlisting}
% Initial guess for Newton iteration
coef=zeros(size(Stiff,1),1);

% Newton iteration
noi=0;
for ii=1:50
noi=noi+1;

  % Compute the nonlinear vector and its Jacobian matrix at
  % the current iterate
  [H,HJ]=nonlinear_global(coef,mesh,@freact_nonlinear,degree);
	
  % Form the residual of the system
  Res = Stiff*coef + H - F;
	
  % Form the Jacobian matrix of the system
  % (w.r.t. unknown coefficients coef)
  J = Stiff + HJ ;
	
  % Solve the linear system for the correction "w"
  w = J \ (-Res);
	
  % Update the iterate
  coef = coef + w;
	
  % Check the accuracy
  if norm(J*w+Res) < 1e-20
    break;
  end
	
end
\end{lstlisting}

\noindent To obtain the non-linear vector $H$ and its Jacobian $HJ$ at the current iterate, we call the function {\it nonlinear\_global}, and it uses the function handle {\it freact\_nonlinear} which is a sub-function in the file {\it Main\_Nonlinear}. The sub-function {\it freact\_nonlinear} has to be supplied by user as the non-linear term $r(u)$ and its derivative $r'(u)$.

\begin{lstlisting}
% Nonlinear reaction
function [r,dr] = freact_nonlinear(u)
  % Value of the nonlinear reaction term at the current iterate
  r = u.^2;
  % Value of the derivative of the nonlinear reaction
  % term at the current iterate
  dr = 2*u;
end
\end{lstlisting}

\section{MATLAB routines for main code}

Here, we give the main m-file {\it Main\_Nonlinear.m} of the code. The full code is available upon request to the e-mail address \email{uzunca@gmail.com}.

\begin{lstlisting}
1 % This routine solves the diffusion-convection-reaction equation
2 %
3 % \alpha u - \epsilon*\Delta u + b\dot\nabla u + r(u) = f
4 %
5 % using DG-FEM.
6
7 function Main_Nonlinear()
8
9 clear all
10 clc
11
12 % Generate the mesh
13
14 % Nodes
15 Nodes = [0,0;0.5,0;1,0;0,0.5;0.5,0.5;1,0.5;0,1;0.5,1;1,1];
16 % Elements
17 Elements = [4,1,5;1,2,5;5,2,6; 2,3,6;7,4,8;4,5,8;8,5,9;5,6,9];
18 % Dirichlet bdry edges
19 Dirichlet = [1,2;2,3;1,4;3,6;4,7;6,9;7,8;8,9];
20 % Neumann bdry edges
21 Neumann = [];
22 % Initial mesh struct
23 mesh = getmesh(Nodes,Elements,Dirichlet,Neumann);
24
25 for jj=1:2
26   mesh=uniformrefine(mesh); %Refine mesh
27 end
28
29 % method : NIPG=1, SIPG=2, IIPG=3
30 method=2;
31
32 % Degree of polynomials
33 degree=1;
34
35 % Set up the problem
36 [penalty,kappa]=set_parameter(method,degree);
37
38 %Compute global matrices and rhs global vector
39 [D,C,R,F]=global_system(mesh,@fdiff,@fadv,@freact,...
40      @fsource,@DBCexact,@NBCexact,penalty,kappa,degree);
41
42 Stiff=D+C+R; % Stiffness matrix
43
44 % Initial guess for Newton iteration
45 coef=zeros(size(Stiff,1),1);
46
47 % Newton iteration
48 noi=0;
49 for ii=1:50
50    noi=noi+1;
51
52    % Compute the nonlinear vector and its Jacobian matrix at
53    % the current iterate
54    [H,HJ]=nonlinear_global(coef,mesh,@freact_nonlinear,degree);
55
56    % Form the residual of the system
57    Res = Stiff*coef + H - F;
58
59    % Form the Jacobian matrix of the system
60    % (w.r.t. unknown coefficients coef)
61    J = Stiff + HJ ;
62
63    % Solve the linear system for the correction "w"
64    w = J \ (-Res);
65
66    % Update the iterate
67    coef = coef + w;
68
69    % Check the accuracy
70    if norm(J*w+Res) < 1e-20
71      break;
72    end
73
74 end
75
76 % Compute L2-error and plot the solution
77 [l2err,hmax]=dg_error(coef,mesh,@fexact,@fdiff,degree);
78
79 % Degree of freedom
80 dof=size(mesh.Elements,1)*(degree+1)*(degree+2)*0.5;
81
82 fprintf(' DoFs h_max L2-error #it\n')
83
84 fprintf('%7d %5.3f %5.3e %d\n',...
85 		dof, hmax ,l2err,noi);
86 end % End of function Main_Nonlinear
87
88 %% Define diffusion, advection, and reaction as subfunctions
89
90 % Diffusion
91 function diff = fdiff(x,y)
92   diff = (10^(-6)).*ones(size(x));
93 end
94
95 % Advection
96 function [adv1,adv2] = fadv(x,y)
97   adv1 =(1/sqrt(5))*ones(size(x));
98   adv2 =(2/sqrt(5))*ones(size(x));
99 end
100
101 % Linear reaction
102 function react = freact(x,y)
103   react = ones(size(x));
104 end
105
106 % Nonlinear reaction
107 function [r,dr] = freact_nonlinear(u)
108   % Value of the nonlinear reaction term at the current iterate
109   r = u.^2;
110
111   % Value of the derivative of the nonlinear reaction
112   % term at the current iterate
113   dr = 2*u;
114 end
115
116 %% Define exact solution and force as subfunctions
117
118 % Exact solution
119 function [yex,yex_x,yex_y] = fexact(fdiff,x,y)
120   % Evaluate the diffusion function
121   diff = feval(fdiff,x,y);
122   % Exact value
123   yex=0.5*(1-tanh((2*x-y-0.25)./(sqrt(5*diff))));
124   % First derivative wrt x
125   yex_x=(-1./(sqrt(5*diff))).*(sech((2*x-y-0.25)./...
126      (sqrt(5*diff)))).^2;
127   % First derivative wrt y
128   yex_y=((0.5)./(sqrt(5*diff))).*(sech((2*x-y-0.25)./...
129      (sqrt(5*diff)))).^2;
130 end
131
132 % Force function
133 function source = fsource(fdiff,fadv,freact,x,y)
134   % Evaluate the diffusion function
135   diff = feval(fdiff,x,y );
136   % Evaluate the advection function
137   [adv1,adv2] = feval(fadv,x, y );
138   % Evaluate the linear reaction function
139   reac = feval(freact,x,y);
140
141   % Exact value
142   yex=0.5*(1-tanh((2*x-y-0.25)./(sqrt(5*diff))));
143   % First derivative wrt x
144   yex_x=(-1./(sqrt(5*diff))).*(sech((2*x-y-0.25)./...
145      (sqrt(5*diff)))).^2;
146   % First derivative wrt y
147   yex_y=((0.5)./(sqrt(5*diff))).*(sech((2*x-y-0.25)./...
148      (sqrt(5*diff)))).^2;
149   % Second derivative wrt x
150   yex_xx=((0.8)./diff).*tanh((2*x-y-0.25)./(sqrt(5*diff))).*...
151      (sech((2*x-y-0.25)./(sqrt(5*diff)))).^2;
152   % Second derivative wrt y
153   yex_yy=((0.2)./diff).*tanh((2*x-y-0.25)./(sqrt(5*diff))).*...
154      (sech((2*x-y-0.25)./(sqrt(5*diff)))).^2;
155
156   % Force function
157   source=-diff.*(yex_xx+yex_yy)+(adv1.*yex_x+adv2.*yex_y)+...
158      reac.*yex+yex.^2;
159 end
160
161
162 %% Boundary Conditions
163
164 % Drichlet Boundary Condition
165 function DBC=DBCexact(fdiff,x,y)
166   % Evaluate the diffusion function
167   diff = feval(fdiff,x,y);
168   % Drichlet Boundary Condition
169   DBC=0.5*(1-tanh((2*x-y-0.25)./(sqrt(5*diff))));
170 end
171
172 % Neumann Boundary Condition
173 function NC = NBCexact(mesh,fdiff,x,y)
174   % Neumann Boundary Condition
175   NC=zeros(size(x));
176 end
177
178
179
180 %% Set-up parameters function for DG-FEM
181
182 function [penalty,kappa]=set_parameter(method,degree)
183
184 global Equation;
185
186 Equation.b0=1; % Superpenalization parameter (In standart b0=1)
187 Equation.base=2; % Choose the basis ( 1:monomials, 2:Dubiner Basis)
188
189 switch method
190   case 1
191     %NIPG
192     Equation.method=1;
193     kappa=1; % type of primal method
194     penalty=1; % penalty parameter
195   case 2
196     %SIPG
197     Equation.method=2;
198     kappa=-1; % type of primal method
199     penalty=3*degree*(degree+1); % penalty parameter
200   case 3
201     %IIPG
202     Equation.method=3;
203     kappa=0; % type of primal method
204     penalty=3*degree*(degree+1); % penalty parameter
205
206 end
207
208 end
\end{lstlisting}

\end{document}